\documentclass[11pt]{amsart}
\usepackage{enumerate,graphicx}

\newcommand{\R}{\mathbb{R}}
\newcommand{\ignore}[1]{}

\newtheorem{theorem}{Theorem}

\newtheorem*{tucker-theorem}{Tucker's Lemma}
\newtheorem*{fan-lemma}{Ky Fan's combinatorial lemma}

\begin{document}

\title{A Constructive Proof of Ky Fan's 
Generalization of Tucker's Lemma}

\author[Timothy Prescott]{Timothy Prescott $^*$}
\thanks{$^*$Research partially supported by a Beckman Research Grant at
Harvey Mudd College.}
\address{Department of Mathematics,
         University of California, Los Angeles, CA 90095}
\email{tmpresco@math.ucla.edu}

\author[Francis Su]{Francis Edward Su $^{**}$}
\thanks{$^{**}$Research partially supported by NSF Grant DMS-0301129.}
\address{Department of Mathematics, Harvey Mudd College,
    Claremont, CA  91711}
\email{su@math.hmc.edu}

\begin{abstract}
We present a constructive proof of Ky Fan's combinatorial lemma concerning 
labellings of triangulated spheres.  Our construction works 
for triangulations of $S^n$ that contain a flag of hemispheres.  
As a consequence, we produce a constructive proof of Tucker's lemma
that holds for a larger class of triangulations than previous constructive
proofs.
\end{abstract}


\maketitle


\section{Introduction}
Tucker's lemma is a combinatorial analogue of the Borsuk-Ulam theorem
with many useful applications.  For instance, it can provide 
elementary routes to proving the Borsuk-Ulam theorem \cite{borsuk} and the 
Lusternik-Schnirelman-Borsuk set covering theorem \cite{freund-todd}, 
Kneser-type coloring theorems \cite{ziegler}, and ``fair division''
theorems in game theory \cite{simmons-su}.  Moreover, any {\em constructive}
proof of Tucker's lemma provides algorithmic interpretations of these results.

Although Tucker's lemma was originally stated for triangulations of an $n$-ball
(for $n=2$ in \cite{tucker} and general $n$ in \cite{lefshetz}), in
this paper 
we shall consider an equivalent version on triangulations of a sphere:
\begin{tucker-theorem}[\cite{tucker,lefshetz}]
Let $K$ be a barycentric subdivision of the octahedral 
subdivision of the $n$-sphere $S^n$.  Suppose that each vertex
of $K$ is assigned a label from $\{\pm 1, \pm 2,... \pm n\}$ in such a
way that labels at antipodal vertices sum to zero.  
Then some pair of adjacent vertices of $K$ have labels that sum to zero.
\end{tucker-theorem}
The original version on the $n$-ball can be obtained from this
by restricting the above statement to a hemisphere of $K$.
This gives a triangulation of the $n$-ball in
which the antipodal condition holds for vertices on the boundary of
the ball.  It is relatively easy to show that Tucker's Lemma is 
equivalent to the Borsuk-Ulam theorem, which says that 
any continuous function $f:S^n\to\R^n$ 
must map some pair of opposite points to the same point 
in the range \cite{borsuk}.  In fact, this equivalence
shows that the triangulation need not be a refinement of
the octahedral subdivision; it need only be symmetric.

However, all known constructive proofs of Tucker's lemma 
seem to require some condition on the triangulation.
For instance, the first constructive proof, 
due to Freund and Todd~\cite{freund-todd},
requires the triangulation to be a refinement of the octahedral
subdivision, and the constructive proof of Yang \cite{yang} depends on 
the $AS$-triangulation that is closely related to the octahedral subdivision.

In this paper, we give a constructive proof of Tucker's lemma for
triangulations with a weaker condition: that it only contain a {\em flag of
hemispheres}.  Our proof (see Theorem \ref{our-tucker})
arises as a consequence of a constructive proof that we develop 
for the following theorem of Fan:

\begin{fan-lemma}[\cite{fan-52}]
Let $K$ be a barycentric subdivision of the octahedral 
subdivision of the $n$-sphere $S^n$.  Suppose that each vertex of 
$K$ is assigned a label from $\{\pm 1, \pm 2,... \pm m\}$ in such a
way that (i) labels at antipodal vertices sum to zero and 
(ii) labels at adjacent vertices do not sum to zero. 
Then there are an odd number of $n$-simplices whose labels are
of the form $\{ k_0, -k_1, k_2, \ldots, (-1)^n k_n \}$, where
$1\le k_0 < k_1 < \dots < k_n \le m$.
In particular, $m\geq n+1$.
\end{fan-lemma}


Our constructive version of Fan's lemma (see Theorem \ref{our-fan}) 
only requires that the triangulation contain a flag of hemispheres.
We use the contrapositive (with $m=n$) 
to obtain a constructive proof of Tucker's lemma.
This yields an algorithm for Tucker's lemma that is
quite different in nature than that of Freund and Todd~\cite{freund-todd}.

Our approach may provide new techniques for developing constructive proofs
of certain generalized Tucker lemmas (such as the $Z_p$-Tucker lemma
of Ziegler \cite{ziegler} or the generalized Tucker's lemma
conjectured by Simmons-Su \cite{simmons-su})
as well as provide new interpretations
of algorithms that depend on Tucker's lemma 
(see \cite{simmons-su} for applications to cake-cutting,
Alon's necklace-splitting problem, 
team-splitting, and other fair division problems).

\ignore{************
We make one comment here on the notion of a constructive proof.  
In a finite setting, one may wonder what is meant by ``constructive'' when
the finite number of possibilities can be checked exhaustively.  By
constructive proof, we mean one that 
(i) shows the existence of the solution and 
(ii) locates it by a method other than an exhaustive search.  
The distinction from an exhaustive search is important for two
reasons.  An exhaustive search is an algorithm, but it cannot
guarantee the existence of a solution without knowing that existence
by some other means.  Secondly, for
continuous results (such as the Borsuk-Ulam theorem) that are obtained
from Tucker's lemma by taking limits as the mesh size of the triangulation 
approaches zero, an exhaustive search is of no help in the limit.  
On the other hand, a constructive proof of
Tucker's lemma for given mesh sizes can be adapted by homotopy methods
to yield algorithms that converge to solutions (such as Borsuk-Ulam 
antipodal points) in a continuous fashion. 
See \cite{todd, yang} for surveys of homotopy methods for
simplicial algorithms.
***************************}

\subsubsection*{Acknowledgements}
The authors are grateful to Joshua Greene for stimulating
conversations related to this work.
%
%

\section{Terminology}
Let $S^n$ denote the $n$-sphere, which we identify with the unit
$n$-sphere $\{x \in \R^{n+1} : \|x\| = 1\}$.
If $A$ is a set in $S^n$, let $-A$ denote the {\em antipodal} set.

A {\em flag of hemispheres} in $S^n$ 
is a sequence $H_0 \subset \dots \subset H_n$
where each $H_d$ is homeomorphic to a $d$-ball, 
and for $1\leq d\leq n$, 
$
\partial H_d = \partial (-H_d) 
= H_d \cap -H_d
= H_{d-1} \cup -H_{d-1}
\cong S^{d-1}
$, 
$H_n\cup-H_n=S^n$, and $\{ H_0,-H_0\}$ are antipodal points.  
One can think of a flag of hemispheres in the
following way: decompose $S^n$ into two balls that intersect
along an equatorial $S^{n-1}$.  Each ball can be thought of as a
hemisphere.  By successively decomposing equators in this fashion
(since they are spheres) and choosing one such ball in each dimension,
we obtain a flag of hemispheres.

A triangulation $K$ of $S^n$ is {\em (centrally) symmetric} if when a simplex
$\sigma$ is in $K$, then $-\sigma$ is in $K$.
A symmetric triangulation of $S^n$ is said to be  
{\em aligned with hemispheres} if we can find a flag of hemispheres 
such that $H_d$ is contained in the $d$-skeleton of the triangulation.
The {\em carrier hemisphere} of a simplex $\sigma$ in $K$ is the 
minimal $H_d$ or $-H_d$ that contains $\sigma$.

A {\em labeling} of the triangulation assigns an integer
to each vertex of the triangulation.
We will say that a symmetric triangulation has an 
{\em anti-symmetric} labeling if each pair of antipodal vertices 
have labels that sum to zero.  We say an edge is a {\em complementary edge}
if the labels at its endpoints sum to zero.

We call a simplex in a labelled triangulation {\em alternating} if 
its vertex labels are distinct in magnitude and 
alternate in sign when arranged 
in order of increasing absolute value, i.e., the labels have the form
$$
\{k_0, -k_1, k_2, \ldots, (-1)^n k_n\} \quad\mbox{ or }\quad 
\{-k_0, k_1, -k_2, \ldots, (-1)^{n+1} k_n\}
$$
where $1 \le k_0 < k_1 < \dots < k_n \le m$.
The {\em sign} of an alternating simplex is the sign
of $k_0$, that is, the sign of the smallest label in absolute value.
For instance, a simplex with labels 
$\{3, -5, -2, 9\}$ is a negative alternating simplex, since the labels
can be reordered $\{-2, 3, -5, 9\}$.  A simplex with labels
$\{-2, 2, -5 \}$ is not alternating because the vertex labels are not
distinct in magnitude.

\begin{figure}[thpb]
\begin{center}
\includegraphics[height=1.25in]{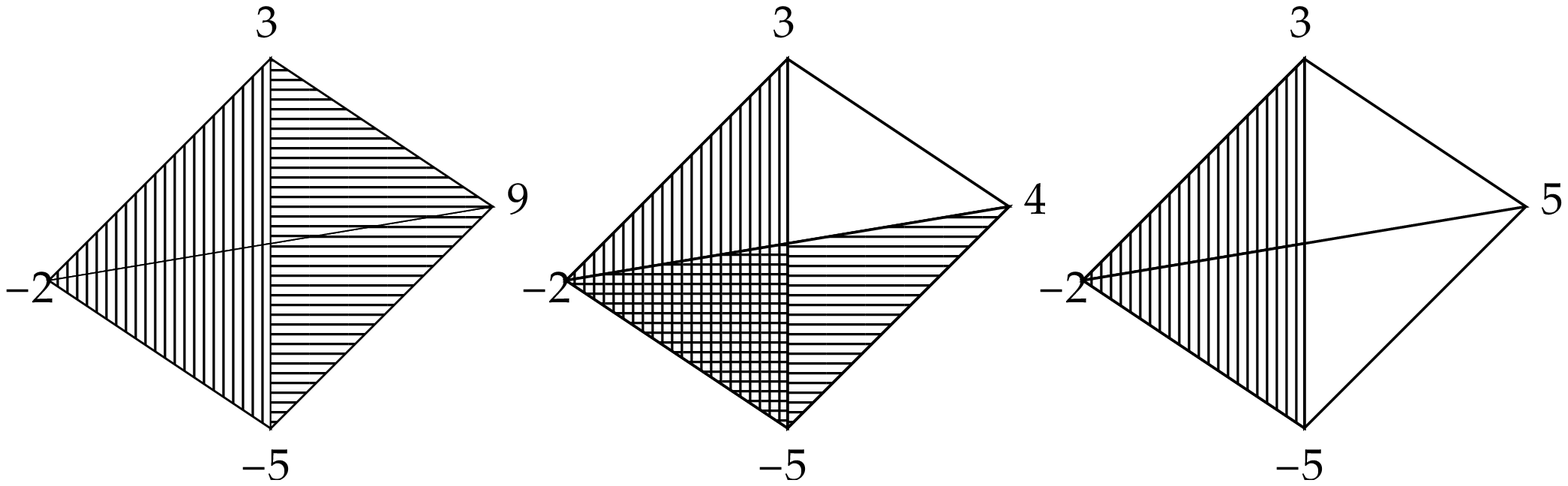}
\end{center}
\caption{The first simplex is alternating and the other two are
  almost-alternating simplices.  Their shaded facets are the
  facets that are also alternating simplices.  The last simplex has
  a complementary edge.}
\label{alternating}
\end{figure}

We also define a simplex to be {\em almost-alternating}
if it is not alternating, but by deleting one
of the vertices, the resulting simplex (a facet) 
is alternating.  The {\em sign} of an almost-alternating simplex
is defined to be the sign any of its alternating facets 
(it is easy to check that this is well-defined).  
For example, a simplex with labels 
$\{-2, 3, 4, -5 \}$ is not alternating, but it is
almost-alternating because deleting $3$ or $4$ would make the
resulting simplex alternating.
Also, a simplex with labels 
$\{-2, 3, 3, -5 \}$ is almost-alternating because deleting
either $3$ would make the resulting simplex alternating.
Finally, a simplex with labels 
$\{-2, 2, 3, -5 \}$ is almost-alternating because deleting 
$2$ would make the resulting simplex alternating.  However, this type of
simplex will not be allowed by the conditions of Fan's lemma 
(since complementary edges are not allowed).  See Figure \ref{alternating}.

Note that in an almost-alternating simplex with no complementary edge,
there are exactly two vertices each of whose removal makes the
resulting simplex alternating, and their labels must be
adjacent to each other when the labels are ordered by
increasing absolute value (e.g., see the second simplex
in Figure \ref{alternating}).  Thus any such almost-alternating simplex must
have exactly two facets which are alternating.

\section{Fan's Combinatorial Lemma}

We now present a constructive proof of Fan's lemma,
stated here for more general triangulations than Fan's original version.

\begin{theorem}
\label{our-fan}
Let $K$ be a symmetric triangulation of $S^n$ aligned with
hemispheres.  Suppose $K$ has
(i) an anti-symmetric labelling by labels $\{\pm 1, \pm 2,... \pm m\}$
and (ii) no complementary edge (an edge whose labels sum to zero).

Then there are an odd number of positive alternating $n$-simplices and
an equal number of negative alternating $n$-simplices.
In particular, $m\geq n+1$.  Moreover, there is a constructive
procedure to locate an alternating simplex of each sign.
\end{theorem}

Fan's proof in \cite{fan-52} used a non-constructive parity argument and 
induction on the dimension $n$.  
Freund and Todd's constructive proof of Tucker's lemma \cite{freund-todd} 
does not appear to generalize to a proof of Fan's lemma, 
since their construction uses $m=n$ in an inherent way.
Cohen~\cite{cohen} implicitly proves a version of Fan's lemma
for $n=2$ and $n=3$ in order to prove Tucker's lemma; his approach 
differs from our proof in that the paths of his search procedure
can pair up alternating simplices with non-alternating simplices
(for instance, $\{1,-2,3\}$ can be paired up with $\{1,-2,-3\})$. 
Cohen hints, but does not explicitly say, 
how his method would extend to higher dimensions; moreover, such an
approach would only be semi-constructive, 
since as he points out, finding one asserted edge in dimension $n$
would require knowing the location of ``all relevant simplices'' 
in dimension $n-1$.

Our strategy for proving Theorem \ref{our-fan} constructively is to identify
paths of simplices whose endpoints are alternating $n$-simplices
or alternating $0$-simplices (namely, $H_0$ or $-H_0$).
Then one can follow such a path from $H_0$ to 
locate an alternating $n$-simplex.

\begin{proof}
Suppose that the given triangulation $K$ of $S^n$ is aligned with the flag 
of hemispheres $H_0 \subset \dots \subset H_n$.
Call an alternating or almost-alternating simplex 
{\em agreeable} if the sign of that simplex matches the sign of its
carrier hemisphere.  For instance, the simplex with labels 
$\{-2, 3, -5, 9\}$ in Figure \ref{alternating}
is agreeable if its carrier hemisphere is $-H_d$ for some $d$.

We now define a graph $G$. 
A simplex $\sigma$ carried by $H_d$ 
is a node of $G$ if it is one of the following:
\begin{enumerate}[(1)]
\item an agreeable alternating $(d-1)$-simplex,
\item an agreeable almost-alternating $d$-simplex, or 
\item an alternating $d$-simplex.
\end{enumerate}

Two nodes $\sigma$ and $\tau$ are adjacent in $G$ if
all the following hold:
\begin{enumerate}[(a)]
\item one is a facet of the other, 
\item $\sigma\cap\tau$ is alternating, and
\item the sign of the carrier hemisphere 
of $\sigma\cup\tau$ matches the sign of $\sigma\cap\tau$.
\end{enumerate}

%

We claim that $G$ is a graph in which every vertex has degree 1 or 2.
Furthermore, a vertex has degree 1 if and only if its 
simplex is carried by $\pm H_0$ or
is an $n$-dimensional alternating simplex.
To see why, we consider the three kinds of nodes in $G$:

\begin{enumerate}[(1)]

\item
An agreeable alternating $(d-1)$-simplex $\sigma$ with carrier $\pm
H_d$ is the facet of exactly two $d$-simplices, 
each of which must be an agreeable alternating or 
an agreeable almost-alternating simplex in the same carrier.
These satisfy the adjacency
conditions (a)-(c) with $\sigma$, hence $\sigma$ has degree 2 in $G$.

\item
An agreeable almost-alternating $d$-simplex $\sigma$ with carrier $\pm H_d$
is adjacent in $G$ to its two facets that are agreeable alternating
$(d-1)$-simplices.  
(Adjacency condition (c) is satisfied because $\sigma$ is agreeable and an
almost-alternating $d$-simplex must have 
the same sign as its alternating facets.)

\item
An alternating $d$-simplex $\sigma$ carried by $\pm H_d$ has one
alternating facet $\tau$ 
whose sign agrees with the sign of the carrier hemisphere
of $\sigma$.  That facet is 
obtained by deleting either the highest or lowest label (by magnitude)
of $\sigma$ so that the remaining simplex satisfies condition (c).
(For instance, the first simplex in Figure \ref{alternating} has two
alternating facets, but only one of them can have a sign that agrees with
the carrier hemisphere.)  Thus $\sigma$ is adjacent to $\tau$ in $G$.

Also, $\sigma$ is the facet of exactly two simplices, 
one in $H_{d+1}$ and one in $-H_{d+1}$, but it is adjacent in $G$ to
exactly one of them; which one 
is determined by the sign of $\sigma$, since the adjacency 
condition (c) must be satisfied.

Thus $\sigma$ has degree 2 in $G$, unless $d=0$ or $d=n$:
if $d=0$, then $\sigma$ is the point $\pm H_0$
and it has no facets, so $\sigma$ has degree 1; 
and if $d=n$, then $\sigma$ is
not the facet of any other simplex, and is therefore of degree 1.
\end{enumerate}

Every node in the graph therefore has degree two 
with the exception of the points at $\pm H_0$
and all alternating $n$-simplices.  Thus $G$ consists of a 
collection of disjoint paths with endpoints at $\pm H_0$ or in
the top dimension.  

Note that the antipode of any path in $G$ is also a path in $G$.
No path can have antipodal endpoints (else the center edge or
node of the path would be antipodal to itself); thus a path is never
identical to its antipodal path.  So all the paths in $G$ must come in
pairs, implying that the number of endpoints of 
paths in $G$ must be a multiple of four.
Since exactly two such endpoints are the nodes at $H_0$ and $-H_0$, 
there are twice an odd number of alternating $n$-simplices.  
And, because every positive alternating $n$-simplex has a
negative alternating $n$-simplex as its antipode, exactly half of the
alternating $n$-simplices are positive.  Thus there are an odd number
of positive alternating $n$-simplices (and an equal number of negative
alternating $n$-simplices).

To locate an alternating simplex, follow the path that begins at
$H_0$; it cannot terminate at $-H_0$ (since a path is never its own 
antipodal path), so it must terminate in a (negative or positive) 
alternating simplex.  The antipode of this simplex will be an
alternating simplex of the opposite sign.
\end{proof}

\begin{figure}[th]
\begin{center}
\includegraphics[height=4in]{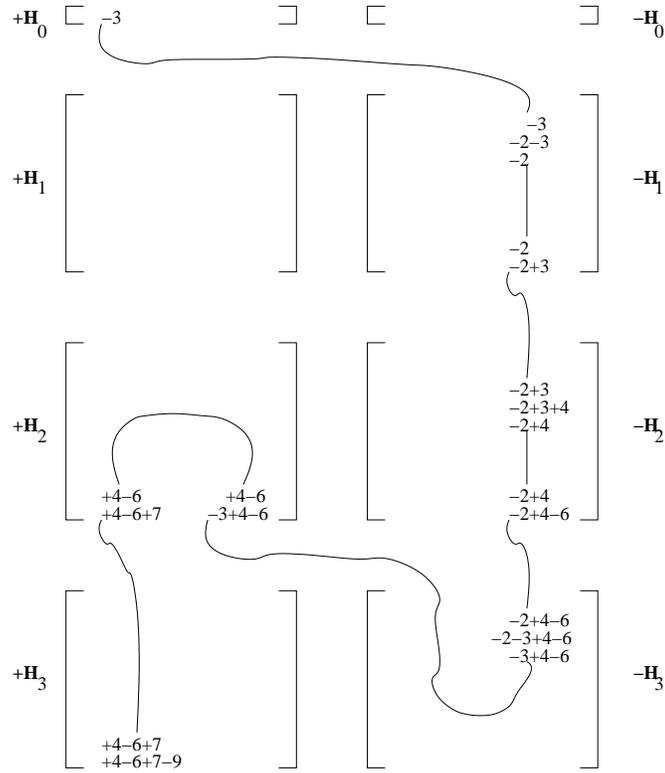}
\end{center}
\caption{An example of what sets of labels of simplices 
            along a path in $G$ could look like.  Repeated labels are
            not shown.}
\label{path}
\end{figure}

Figure \ref{path} shows an example of how a path
may wind through the various hemispheres of a triangulated $3$-sphere.
Note how the sign of each simplex agrees with the sign of its
carrier hemisphere (``agreeability''), 
unless the path connects a $d$-hemisphere with a
$d+1$-hemisphere, in which case the sign of the $d$-simplex specifies
which $(d+1)$-hemisphere the path should connect to.  
These facts follow from adjacency condition (c).

Our approach is related to that of another paper of Fan \cite{fan-67}, 
which studied labelled triangulations of an $d$-manifold $M$ 
and derived a set of paths that pair up alternating simplices 
in the interior of $M$ with positive alternating simplices on the 
boundary of $M$.  When $M=H_d$, the paths of Fan
coincide with the restriction of our paths in $G$ to $H_d$.  
By itself, this is only semi-constructive, 
since finding one alternating $d$-simplex necessitates locating all 
positive alternating $(d-1)$-simplices on the boundary of $H_d$.  
To make Fan's approach fully constructive for $S^n$, 
one might attempt to 
use Fan's approach in each $d$-hemisphere of $S^n$ and then
glue all the hemispheres in each dimension together, thereby gluing
all the paths.  But this results in paths that branch (where positive
alternating simplices in $H_d$ are glued to paths in {\em both} $H_{d+1}$
and $-H_{d+1}$) or paths that terminate prematurely (where a path ends
in a negative alternating $d$-simplex where $d<n$).

By contrast, the path we follow in $G$ from the point $H_0$ to an
alternating $n$-simplex is well-defined, 
has no branching, and need not pass through 
all the alternating $(d-1)$-simplices on the boundary of $H_d$ for
each $d$.
In our proof, the use of the flag of hemispheres
controls the branching that would occur in paths of $G$ if 
one ignored the property of ``agreeability'' and adjacency condition
(c).  In that sense, it serves a similar function in controlling
branching as the 
use of the flag of polytope faces in the constructive proof of the polytopal
Sperner lemma of DeLeora-Peterson-Su \cite{deloera-peterson-su}.

Note that the 
contrapositive of Theorem \ref{our-fan} implies Tucker's
Lemma, since if $m=n$ and condition (i) holds, then condition (ii)
must fail.  In fact, if we remove condition (ii) in the statement of Theorem
\ref{our-fan}, the graph $G$ can
have additional nodes of degree 1, namely, agreeable almost-alternating
simplices with a complementary edge.

This gives a constructive
proof for Tucker's lemma by starting at $H_0$ and following the
associated path in $G$. Because $m=n$, there are not enough labels for
the existence of any alternating $n$-simplices, so there must be an
odd number of agreeable positive almost-alternating simplices with a
complementary edge.  (Note that this says nothing 
about the parity of the number of complementary edges,
since several such simplices could share one edge.)

It is of some interest that our constructive proof allows for a larger
class of triangulations than previous constructive 
proofs of Tucker's lemma, so for completeness we state it carefully here:

\begin{theorem}
\label{our-tucker}
Let $K$ be a symmetric triangulation of $S^n$ aligned with
hemispheres.  Suppose $K$ has an anti-symmetric labelling by 
labels $\{\pm 1, \pm 2, \ldots, \pm n\}$.  Then there are an odd
number of positive (negative) almost-alternating simplices which
contain a complementary edge.  Moreover, there is a constructive
procedure to locate one such edge.
\end{theorem}

The hypothesis that $K$ can be aligned with a flag of hemispheres is
weaker than, for instance, requiring $K$ to refine the octahedral
subdivision (e.g., Freund-Todd's proof of Tucker's lemma).
If a triangulation refines the octahedral subdivision,
then the octahedral orthant hyperplanes contain a natural flag of
hemispheres.  But the converse is not true: 
there are triangulations aligned with 
hemispheres that are not refinements of the octahedral subdivision.
For instance, consider 
the triangulated $2$-sphere $\{(x,y,z):x^2+y^2+z^2=1\}$
whose $1$-skeleton is cut out by intersections with 
the plane $z=0$ and half-planes
$\{x=0, z\geq 0\}$ and 
$\{y=0, z\leq 0\}$.  This triangulation has 4 vertices at $(\pm 1,0,0)$
and $(0,\pm 1,0)$, it is symmetric, and it contains a flag of hemispheres.
But it is combinatorially equivalent to the boundary of a $3$-simplex, 
and hence does not refine the octahedral subdivision of $S^2$.

We remark that the $AS$-triangulation, used by Yang \cite{yang} to
prove Tucker's lemma, is closely related to an octahedral subdivision
and contains a natural flag of hemispheres.

It is an interesting open question as to whether any symmetric
triangulation of $S^n$ can be aligned with a flag of hemispheres, and
if so, how to find such a flag.  Together with our arguments this
would yield a constructive proof of Tucker's lemma for any symmetric 
triangulation.

\bibliographystyle{abbrv}
\bibliography{fixedpoint}

\end{document}